\documentclass{article}

\usepackage{amsthm}
\usepackage{amsfonts}
\usepackage{amsmath}
\usepackage{amssymb}
\usepackage{latexsym}
\usepackage{graphicx}

 \theoremstyle{plain}
 \newtheorem{thrm}{Theorem}[section]

 \theoremstyle{definition}

 \theoremstyle{remark}
 
\newcommand{\comment}[1]{}

\begin{document}

\title{A note on
Beckett-Gray codes and the relationship of Gray codes to data structures}

\author{Mark Cooke \\
\small Network Appliance \\[-0.8ex]
\small 495 East Java Drive \\[-0.8ex]
\small Sunnyvale, CA 94089. U.S.A.\\[-0.8ex]
\small \texttt{Mark.Cooke@netapp.com}
\and
Chris North \\ 
\small Tutte Institute for Mathematics and Computer Science \\[-0.8ex]
\small 1929 Ogilvie Road \\ [-0.8ex]
\small Ottawa ON K1J 0B9, Canada \\[-0.8ex]
\small \texttt{cjnorth@gmail.com}
\and
 Megan Dewar, Brett Stevens\\
\small School of Mathematics and Statistics \\[-0.8ex]
\small Carleton University \\ [-0.8ex]
\small 1125 Colonel By Drive \\ [-0.8ex]
\small Ottawa ON K1S 5B6, Canada  \\[-0.8ex]
\small \texttt{megan.dewar@gmail.com}\\[-0.8ex]
\small \texttt{brett@math.carleton.ca}}
\maketitle

\begin{abstract} 
In this paper we introduce a natural mathematical structure derived from Samuel Beckett's play ``Quad''.  We call this structure a binary Beckett-Gray code.  We enumerate all codes for $n \leq 6$ and give examples for $n=7,8$. Beckett-Gray codes can be realized as successive states of a queue data structure.  We show that the binary reflected Gray code can be realized as successive states of two stack data structures.
\end{abstract}

\section{Introduction}

Samuel Beckett's play ``Quad'' consists of a series of arrivals and departures of four characters resulting in their appearance on stage together in different combinations throughout the play.  At regular intervals exactly one character will enter the stage or exactly one character will exit the stage.  There are no apparent constraints on which off-stage character is allowed to enter, however, when a character exits, it is always the character who has been on stage the longest.  The play ends as one character's exit is about to reproduce the play's starting configuration. Beckett's text explicitly notes that every possible non-empty subset of characters appear on stage together at least once but admits that these subsets do not appear a uniform number of times \cite{plays}.

If each subset of appeared precisely once, then the series of subsets -- represented by length four binary words -- would form a cyclic Gray code with $n=4$.  Gray codes of binary $n$-bit words are equivalent to Hamilton paths and cycles in the $n$-dimensional hypercube \cite{MR993775}.  The additional exit constraint is equivalent to only permitting a 1 in bit position $p$ to change to a 0 if, of those positions containing a 1, position $p$ has been so longest. In other words, this 1 has the longest current {\em run}. We will call a Gray code with this additional property a {\em Beckett-Gray code}. For all but the shortest codes, we will display a Gray code by its {\em transition sequence}; the list of bit positions at which each change takes place. We number the bit positions from the right, starting from zero.  While the literary use of Beckett-Gray codes has been explored \cite{brett:beck4}, here we investigate these objects mathematically. We will discuss their isomorphisms and enumerate all Beckett-Gray codes for $n \leq 6$ and give examples for $n=7,8$. We finish by discussing the compatibility of Gray codes -- not just Beckett-Gray codes -- with various data-structures.  This work was completed over a decade ago and since then other results regarding Beckett-Gray codes have been published.  We feel this note provides the necessary background to this work; for a more detailed account of the history and searches see \cite{2016arXiv160806001v1C}.

\section{Preliminaries} \label{mathprelim}

The existence of a Beckett-Gray code is equivalent to  being able to realize every subset of an $n$-set exactly once as successive states of a queue, where queue entries are the bit positions that are currently 1. In Table~\ref{n3} we give an example of a non-cyclic Beckett-Gray code and the corresponding queue states.


\begin{table}[htp]
\begin{center}
\begin{tabular}{|c|l|} 
\hline
\bf{Code}	& \bf{Queue}	\\ \hline \hline
$000$ &	$\emptyset$ \\
$001$ &	$0$ \\
$011$ &	$0,1$ \\
$010$ &	$1$ \\
$110$ &	$1,2$ \\
$100$ & $2$ \\
$101$ &	$2,0$ \\
$111$ &	$2,0,1$ \\
\hline
\end{tabular}
\end{center}
\caption{A 3-bit non-cyclic  binary Beckett-Gray code. \label{n3}}
\end{table}

Two Gray codes are called {\em isomorphic} if one can be obtained from the other by permuting bit positions, adding a fixed $n$-bit word to all words and/or reversing the list of words.  Permuting bit positions preserves the Beckett-Gray property, but addition of a non-zero word switches the role of 0 and 1 in at least one position and may destroy the Beckett-Gray property. On the other hand, given that the time reversal of a queue is still a queue, the reverse of a Beckett-Gray code is also Beckett-Gray.

A Gray code is called {\em self-isomorphic} if there is a non-trivial isomorphism that fixes the code.  A Gray code is {\em self-reverse} if it can be obtained from its reversal by any isomorphism.  We show that no binary Gray code can be isomorphic to its own reverse without adding a fixed binary word.
\begin{thrm} \label{self-rev}
No cyclic binary Gray code is self-reverse for $n \geq 3$ without adding a fixed word.
\end{thrm}
\begin{proof}
Assume we have a such a self-reverse code.  Since the reversal is isomorphic, there exists a permutation of the bit positions, $\rho$, which is the isomorphism between the code and its reverse. Note that $\rho$ is order two and all its orbits are of size 1 or 2.  Let $O \subseteq \{1,2,\ldots,n\}$ be any non-empty union of orbits of $\rho$ and let $w_{O}$ be the word that has a 1 only in bits from $O$. Since $\rho(w_O) = w_{O}$, $w_O$ must occur at position $2^{n-1}$ in the code but this contradicts the fact that when $n \geq 3$, there is more than one choice for $O$.
\end{proof}
If the addition of an $n$-bit binary word is permitted as an isomorphism, then Frank Gray's binary reflected Gray code is self-reverse for all $n$.

In the course of enumerating Beckett-Gray codes, we must be careful to eliminate isomorphic codes.  All codes listed in this paper will be given as the lexicographically least transition sequence from their isomorphism class.

 Donald Knuth included a description of Beckett-Gray codes in Pre-Fascicle 2a of his Vol 4 of {\em The Art of Computer Programming} and posed the problem of finding a Beckett-Gray code for $n=8$ \cite{knuth:09}.  This note presents the 2005 solution to this question.  In 2007, Dennis Wong and Joe Sawada independently enumerated the Beckett-Gray codes for $n=6$ using faster techniques and produced almost 10,000 cyclic Beckett-Gray codes for $n=7$ \cite{sawada:07}.  These results were part of Wong's M.Sc thesis \cite{wong07}.

\section{Enumeration and existence} \label{enum_exist}

The unique binary Gray codes for 1 and 2-bits are both cyclic and are both Beckett-Gray codes.  There are no cyclic Beckett-Gray codes for $n=3$ nor $n=4$, but there are one and four non-cyclic Beckett-Gray codes, respectively.  Their transition sequences are 0102101, 010213202313020, 010213212031321, 012301202301230, and 012301213210321. For $n=5$ there are eight non-isomorphic cyclic Beckett-Gray codes. These codes are given in Table \ref{eg_5}.
\begin{table}[ht!]
\begin{center}
\begin{tabular}{cc}
01020132010432104342132340412304 &
01020312403024041232414013234013 \\
01020314203024041234214103234103 &
01020314203240421034214130324103 \\
01020341202343142320143201043104 &
01023412032403041230341012340124 \\
01201321402314340232134021431041 &
01203041230314043210403202413241 
\end{tabular}
\caption{Transition sequences of the eight cyclic Beckett-Gray codes for $n=5$.}
\label{eg_5}
\end{center}
\end{table}
There are 116 non-isomorphic non-cyclic 5-bit Beckett-Gray codes, available upon request. For $n=6$, there are $94,841$ non-isomorphic cyclic Beckett-Gray codes.  To give some idea of the range of codes, the lexicographically first and last are
\begin{align*}
  &0102013120240312152430145052341304513534523514302514523405125415, \\
  &0123450123435432543125340134140503214541052401432501435032125032
\end{align*}
There are $5,868,331$ non-isomorphic non-cyclic 6-bit Beckett-Gray codes.

There exist cyclic Beckett-Gray codes for $n=7$ and $n=8$.  For $n=7$, the lexicographically first cyclic Beckett-Gray code is
\[\begin{split}
012&3450106212343540616521234640561324201560323415051306451210626 \\
&4215606432150563641505646052563625410363410502350412342104320545.
\end{split}\]
A cyclic Beckett-Gray code for $n=8$ is
\[\begin{split}
012&3456070121324356576071021353462670153741236256701731426206570 \\
&1342146560573102464537571020435376140736304642737035640271327505 \\
&4121027564150240365425013602541615604312576032572043157624321760 \\
&4520417516354767035647570625437242132624161523417514367143164314.
\end{split}\]

The enumeration was done using a depth first lexicographical search algorithm.  To obtain all codes for $n=6$, the algorithm was parallelized and distributed over a small cluster of 2003 era computers. This search was confirmed by being written two times independently and has also been verified by Wong and Sawada \cite{sawada:07}. We used Knuth's method \cite{MR0373371} to determine an estimate for the size of the full back-tracking search tree for $n=7$.  The method was verified for accuracy on the $n=5,6$ cases.  We conclude that the sizes of the trees for $n=5,6$ and $7$ are on the order of $2^{19}$, $2^{44.7}$ and $2^{102}$ respectively. We suspect that an exhaustive tree-based search for the $n=7$ cyclic Beckett-Gray codes, even with advances in computer hardware and the techniques of Sawada and Wong \cite{sawada:07}, is infeasible.  The method used to find cyclic Beckett-Gray codes for $n=8$ used simulated annealing to produce long, incomplete codes that satisfied the Beckett-Gray property.  These were then fed as seeds to a deterministic backtracking search.  This is similar to the approach used in \cite{MR1377601,MR1935751}.

\section{Gray codes and data structures}

The connection between Beckett-Gray codes and queues prompts the question of whether other Gray codes can be realized as successive states of various data structures.  We have been able to show a nice result regarding Frank Gray's binary reflected Gray code \cite{MR0094273,gray:53}.
\begin{thrm}
\label{brg_stack}
The standard binary reflected Gray code can be realized as successive states of a pair of stacks.
\end{thrm}
\begin{proof}
This can be proved by induction.  The two stacks will be used to represent the appearances of a 1 in
even and odd bit positions, respectively.  We show this for small $n$ in Table~\ref{2stack}.  
\begin{table}[ht!]
\begin{tabular}{|l|l|l|c|l|l|l|}
\cline{1-3} \cline{5-7}
 Gray code & even stack & odd stack  && Gray code & even stack & odd stack \\
\cline{1-3} \cline{5-7} 
00 & $\emptyset$ & $\emptyset$  && 000 & $\emptyset$ & $\emptyset$ \\
01 & 0 & $\emptyset$ && 001 & 0 & $\emptyset$ \\
11 & 0 & 1 && 011 & 0 & 1 \\
10 & $\emptyset$ & 1  && 010 & $\emptyset$ & 1 \\
\cline{1-3} 
\multicolumn{3}{c}{} && 110 & 2 & 1 \\
\multicolumn{3}{c}{} && 111 & 20 & 1 \\
\multicolumn{3}{c}{} && 101 & 20 & $\emptyset$ \\
\multicolumn{3}{c}{} && 100 & 2  & $\emptyset$ \\
\cline{5-7}
\end{tabular}
\caption{Binary reflected Gray codes realized as successive states of two stacks.}
\label{2stack}
\end{table}  
Assume that $n=k$ can be realized in two
stacks, the first using even indices and the second using odd indices.  The first $2^k$
words of the $(k+1)$-bit binary reflected Gray code are the same words as
those of the $k$-bit binary reflected Gray code with a 0 in the  $(k+1)^{st}$ bit.  Thus the first half of the code is realizable
in two stacks by the induction hypothesis.  The last word in this half of the code
is all 0s except there is a 1 in bit position $k-1$.  This corresponds to
the stack with parity $k-1$ containing just the element $k-1$ and the
other stack being empty.  The remaining half of the $(k+1)$-bit binary reflected Gray code is
the first half in reverse order, with a 1 always in bit
position $k$.  We simply push $k$ onto the stack with parity $k$ (which is empty) and then notice
that we can easily run the stack operations from the first half in
reverse as stack operations are time reversible.
\end{proof}
We note that the stack data structure has been used in efficient generation of many Gray codes, however,  in these cases the state of the stack itself does not represent the current element of the code \cite{MR0424386}.

\section{Conclusion} \label{conc}

Beckett-Gray codes have connections to other combinatorial orderings.  For example, they force the runs of 1s in the code to be relatively long and so are related to some known Gray codes with long bit runs \cite{goddyn:88,MR2014514}.  Since the objects in a Beckett-Gray code are realizable as the successive states of a queue, these codes are similar to de Bruijn cycles and universal cycles \cite{diaconis:92}, except that in the case of Beckett-Gray codes, successive objects are close both in the Hamming distance sense and in the queue state sense.

Small Beckett-Gray codes all exist, suggesting that the larger codes do too. The goal is to find a recursive or direct construction for Beckett-Gray codes; however, a non-constructive or probabilistic proof would be encouraging.

\bibliographystyle{hplain}
\bibliography{journals,bibliography,lit}

\end{document}